\documentclass[12pt]{article}                                               

\input{amssym.def}                                                           
\input{amssym.tex}                                                           
                                                                             
\begin{document}                                                             
\title{Foliations on modular curves}

\author{Igor  ~V. ~Nikolaev
\footnote{Partially supported 
by NSERC.}}


\date{}
 \maketitle


\newtheorem{thm}{Theorem}
\newtheorem{lem}{Lemma}
\newtheorem{dfn}{Definition}
\newtheorem{rmk}{Remark}
\newtheorem{cor}{Corollary}
\newtheorem{cnj}{Conjecture}
\newtheorem{exm}{Example}


\newcommand{\N}{{\Bbb N}}
\newcommand{\mod}{\hbox{\bf mod}}
\newcommand{\GCD}{\hbox{\bf GCD}}
\newcommand{\GL}{\hbox{\bf GL}}
\newcommand{\SL}{\hbox{\bf SL}}
\newcommand{\rank}{\hbox{\bf rank}}
\newcommand{\End}{\hbox{\bf End}}
\newcommand{\Per}{\hbox{\bf Per}}
\newcommand{\ch}{\hbox{\bf char}}


\centerline{{\it  In memory of  Marco ~Brunella}}

\bigskip

\begin{abstract}
It is proved, that a foliation on a modular curve  given by 
 the vertical trajectories of holomorphic differential corresponding to   the  Hecke eigenform
is either the Strebel foliation or the pseudo-Anosov foliation.

\vspace{7mm}

{\it Key words and phrases:  Hecke eigenforms,  singular foliations}

\vspace{5mm}
{\it MSC: 11F12 (automorphic forms);  57R30 (foliations)}
\end{abstract}

\section{Introduction}
Let $N>1$ be a natural number and  consider a finite index subgroup 
of the modular group given by the formula: 
\begin{equation}
\Gamma_0(N) = \left\{\left(\matrix{a & b\cr c & d}\right)\in SL(2,{\Bbb Z})~|~
c\equiv 0~mod ~N\right\}.
\end{equation}
Let ${\Bbb H}=\{z=x+iy\in {\Bbb C} ~|~ y>0\}$ be the upper half-plane  and 
let $\Gamma_0(N)$  act on ${\Bbb H}$  by   linear fractional
transformations;  consider the  orbifold  ${\Bbb H}/\Gamma_0(N)$.
To compactify the orbifold 
at the cusps, one adds a boundary to ${\Bbb H}$,  so that 
${\Bbb H}^*={\Bbb H}\cup {\Bbb Q}\cup\{\infty\}$ and the compact Riemann surface 
$X_0(N)={\Bbb H}^*/\Gamma_0(N)$ is called a {\it modular curve}.   
The holomorphic functions $f(z)$ on ${\Bbb H}$ that
vanish at the cusps and such that
\begin{equation}
f\left({az+b\over cz+d}\right)=(cz+d)^{2}f(z),\qquad
\forall \left(\matrix{a & b\cr c & d}\right)\in\Gamma_0(N), 
\end{equation}
are  called  {\it cusp forms} of weight two;  the (complex linear) space of such forms
will be denoted by $S_2(\Gamma_0(N))$.  The formula $f(z)\mapsto \omega=f(z)dz$ 
defines an isomorphism  $S_2(\Gamma_0(N))\cong \Omega_{hol}(X_0(N))$, where 
$\Omega_{hol}(X_0(N))$ is the space of holomorphic differentials
on the Riemann surface $X_0(N)$.  Note that 
$\dim_{\Bbb C}(S_2(\Gamma_0(N))=\dim_{\Bbb C}(\Omega_{hol}(X_0(N))=g$,
where $g=g(N)$ is the genus of the surface $X_0(N)$. 
A Hecke operator, $T_n$, acts on $S_2(\Gamma_0(N))$ by the formula
$T_n f=\sum_{m\in {\Bbb Z}}\gamma(m)q^m$, where
$\gamma(m)= \sum_{a| GCD(m,n)}a c_{mn/a^2}$ and 
$f(z)=\sum_{m\in {\Bbb Z}}c(m)q^m$ is the Fourier
series of the cusp form $f$ at $q=e^{2\pi iz}$.  For $(n,N)=1$ it is known that $T_n$ is a
self-adjoint linear operator on the vector space $S_2(\Gamma_0(N))$
endowed with the Petersson inner product;  the algebra
${\Bbb T}_N :={\Bbb Z}[T_1,T_2,\dots]$ is a commutative algebra.
Any cusp form $f_N\in S_2(\Gamma_0(N))$ that is an eigenvector
for one   of $T_n$ is referred  to
as a {\it Hecke eigenform}; such an  eigenform is called
{\it rational}  whenever its Fourier coefficients $c(m)\in {\Bbb Z}$ after the normalization $c(1)=1$.

Let $\phi=\Re~(\omega)$ be the real part
of $\omega$; it is a closed form on the surface $X_0(N)$.
(Alternatively, one can take for $\phi$ the  imaginary part of $\omega$.)
Note that $\phi$ is well-defined modulo a real multiple $C$;
indeed,  because $f_N(z)=\sum_{m\in {\Bbb Z}}c(m)q^m$
and $c(m)\in {\Bbb R}$ it is easy to deduce that for each eigenform
$Cf_N(z)$ the constant $C\in {\Bbb R}$.   
By ${\cal F}$ we understand a measured foliation given by  trajectories
of the  closed  form $\phi$;  we address the following

\medskip\noindent
{\sf Classification problem.} 
{\it Find the topological types of  foliation ${\cal F}_N$
given by the  Hecke eigenform $f_N$.}

\medskip\noindent
Recall that measured foliation ${\cal F}$ 
is called a {\it Strebel foliation}, if all but a finite number of its leaves are closed [Strebel 1984]  \cite{S};
the Strebel foliation is invariant of a parabolic automorphism of the surface. 
The measured foliation ${\cal F}$ is called a {\it pseudo-Anosov foliation}, if it 
is minimal and  invariant  of a hyperbolic (pseudo-Anosov) automorphism of the
surface [Casson \& Bleiler 1988]  \cite{CB}; our main result  can be expressed as follows.
\begin{thm}\label{thm1}
The foliation ${\cal F}_N$ is topologically conjugate to a pseudo-Anosov
foliation unless $f_N$ is a rational eigenform,  in which case it is topologically
conjugate to a  Strebel foliation. 
\end{thm}
The article is organized as follows. All preliminary results are introduced in Section 2;
all  proofs are given in Section 3.

\section{Preliminaries}
\subsection{Measured foliations}
By a $p$-dimensional,  class $C^r$ foliation of an $m$-dimensional 
manifold $M$ one understands a decomposition of $M$ into a union of
disjoint connected subsets $\{ {\cal L}_{\alpha}\}_{\alpha\in A}$, called
the {\it leaves} of the foliation. The leaves must satisfy the
following property: Every point in $M$ has a neighborhood $U$
and a system of local class $C^r$ coordinates 
$x=(x^1,\dots, x^m): U\to {\Bbb R}^m$ such that for each leaf 
${\cal L}_{\alpha}$, the components of $U\cap {\cal L}_{\alpha}$
are described by the equations $x^{p+1}=Const, \dots, x^m=Const$.
Such a foliation is denoted by ${\cal F}=\{ {\cal L_{\alpha}}\}_{\alpha\in A}$.
The number $q=m-p$ is called a {\it codimension} of the foliation
${\cal F}$,  see [Lawson 1974]  \cite{Law1} p.370.    
The codimension $q$ $C^r$-foliations ${\cal F}_0$ and ${\cal F}_1$ 
are said to be $C^s$-conjugate ($0\le s\le r$) if there exists a diffeomorphism
of $M$, of class $C^s$, which maps the leaves of ${\cal F}_0$
onto the leaves of ${\cal F}_1$. If $s=0$, ${\cal F}_0$ and ${\cal F}_1$
are {\it topologically conjugate}   [Lawson 1974]  \cite{Law1},  p.388;
we shall write such foliations as  ${\cal F}_0\cong {\cal F}_1$.   
The foliation ${\cal F}$ is called {\it singular} if the codimension
$q$ of the foliation depends on the leaf.  We further  assume
that $q$ is constant for all  but a {\it finite}
number of leaves. Such a set of the exceptional leaves will be denoted by
$Sing ~{\cal F}:= \{ {\cal L}_{\alpha}\}_{\alpha\in F}$, where $|F|<\infty$.
Note that in the case $F$ is an empty set, one gets the usual definition of a 
(non-singular) foliation. A quick  example of the singular foliations is given
by the trajectories of a non-trivial differential form  on the manifold $M$,
which vanish in a finite number of  points of $M$; 
the set of zeroes of such a form corresponds to the  exceptional leaves of 
the  foliation.      
A {\it measured foliation} ${\cal F}$ on the surface $X$
is the partition of $X$ into the singular points $x_1,\dots,x_n$ of
order $k_1,\dots, k_n$ and regular leaves (1-dimensional submanifolds). 
On each  open cover $U_i$ of $X-\{x_1,\dots,x_n\}$ there exists a non-vanishing
real-valued closed 1-form $\phi_i$  such that: 
(i)  $\phi_i=\pm \phi_j$ on $U_i\cap U_j$;
(ii) at each $x_i$ there exists a local chart $(u,v):V\to {\Bbb R}^2$
such that for $z=u+iv$, we have  $\phi_i=Im~(z^{k_i\over 2}dz)$ on
$V\cap U_i$ for some branch of $z^{k_i\over 2}$.
The pair $(U_i,\phi_i)$ is called an atlas for the measured foliation ${\cal F}$.
Finally, a measure $\mu$ is assigned to each segment $(t_0,t)\in U_i$, which is  transverse to
the leaves of ${\cal F}$, via the integral $\mu(t_0,t)=\int_{t_0}^t\phi_i$. The 
measure is invariant along the leaves of ${\cal F}$, hence the name. 
Note that the measured foliation can have  singular points  with an odd number 
of the saddle sections (half-integer index);  those cannot be continuously  oriented along the leaves
and therefore cannot be given by the trajectories of a closed form $\phi$.
However, each  measured  
foliation can be rendered oriented by taking if necessary  the double cover of $X$ ramified 
at the singular points of the half-integer index; further we work with the foliations given by trajectories
of a closed form $\phi$. 
Let $\{\gamma_1,\dots,\gamma_n\}$ be a
basis in the relative homology group $H_1(X, Sing~{\cal F}; {\Bbb Z})$,
where $Sing~{\cal F}$ is the set of singular points of the foliation ${\cal F}$; 
here  $n=2g+m-1$, where $g$ is the genus of $X$ and $m= |Sing~({\cal F})|$.
(In what follows it will be shown that $m=1$ and an involution on  $H_1(X, Sing~{\cal F}; {\Bbb Z})$
will imply that $n=g$, see proof of lemma \ref{lm6}.)
Let $\lambda_i$ be the periods of $\phi$  against the basis $\gamma_i$.
By a {\it Jacobian}  of  measured foliation ${\cal F}$  we
understand a ${\Bbb Z}$-module $Jac~({\cal F}):={\Bbb Z}\lambda_1+\dots+{\Bbb Z}\lambda_n$
regarded as a subset of the real line ${\Bbb R}$.  
The following lemma describes basic properties  of the Jacobians;
it will be proved in Section 3.1.   
\begin{lem}\label{lm1}
For each uniquely ergodic measured  foliation  ${\cal F}$  the  following is true:

\medskip
 (i)  $Jac~({\cal F})$ is a well-defined ${\Bbb Z}$-module of the real line ${\Bbb R}$;

\smallskip
(ii) if ${\cal F}\cong {\cal F}'$ are topologically conjugate foliations, 
then  $Jac~({\cal F}')=\mu ~Jac~({\cal F})$ for a constant $\mu>0$,
i.e. the projective class of $Jac~({\cal F})$ is an invariant of
topological conjugation. 
\end{lem}

\subsection{Automorphisms of surfaces}
Let $\varphi: X\to X$ be an orientation-preserving automorphism of a compact surface $X$; 
the composition of (the isotopy classes of) such automorphisms define a countable
group denoted by $Mod~(X)$. The measured foliations play an outstanding r\^ole in
the classification of surface automorphisms;  to give an idea, let $X=T^2$
be the two-dimensional torus.  
One can regard $T^2$ as the quotient of the Euclidean plane 
${\Bbb R}^2$ by the integer lattice ${\Bbb Z}^2$, endowed with a fixed
orientation. Since $\pi_1(T^2)={\Bbb Z}\oplus {\Bbb Z}$, the automorphisms
of $T^2$ correspond to the elements of group $GL_2({\Bbb Z})$ as any
element $A\in GL_2({\Bbb Z})$ maps ${\Bbb Z}^2$ to itself and so
induces a continuous map $\varphi_A: T^2\to T^2$;
the map $\varphi_A$ has an inverse $\varphi_A^{-1}$ given by the matrix $A^{-1}$. 
The map $\varphi_A$ preserves orientation if and only if $det~(A)=1$; 
thus $Mod~(T^2)\cong SL_2({\Bbb Z})$.  
If $A\in SL_2({\Bbb Z})$, the characteristic  polynomial of $A$
can be written as $\lambda^2-(tr~(A))\lambda+1$. The eigenvalues 
are either: (i) complex, i.e. $tr~(A)=0$ or $\pm 1$;
(ii) double reals $+1$ or $-1$, i.e. $tr~(A)=2$ or $-2$;
(iii) distinct reals, i.e. $|tr~(A)|>2$.
 At the level of surface automorphisms, the above cases
correspond to: (i) a finite order automorphism $(\varphi_A)^{12}=1$; 
(ii)  matrix $A$ has an integral eigenvector which defines a measured
foliation ${\cal F}$ on $T^2$.  The leaves of ${\cal F}$ are closed and
invariant under the map $\varphi_A$; the latter  is a Dehn twist in $C\in {\cal F}$.
 The automorphism $\varphi_A$ is parabolic  and ${\cal F}$  
is called a {\it Strebel foliation}.
(iii)  $\lambda_A>1>1/\lambda_A$, where $\lambda_A$
is the Perron-Frobenius eigenvalue and we let $v_A, v_A'$ be the
corresponding eigenvectors. Since $\theta_A=v_A^{(2)}/v_A^{(1)}$
is irrational, $\varphi_A$ leaves no simple closed curve invariant.
The translation of vectors $v_A$ and $v_A'$ yields two measured foliations
${\cal F}_u$ and ${\cal F}_s$, which are carried by $\varphi_A$ to the
foliations $\lambda_A{\cal F}_u$ and ${1\over\lambda_A}{\cal F}_s$,
respectively. In other words, $\varphi_A$ is a linear homeomorphism which
stretches by a factor $\lambda_A$ in one direction and shrinks 
by the same factor in a complementary direction. The homeomorphism
$\varphi_A$ is hyperbolic and invariant foliations ${\cal F}_s$ and ${\cal F}_u$ are
 called {\it Anosov foliations}.   What is remarkable, the classification
holds for the higher genus surfaces;  due to natural restrictions we can only outline these results referring
the interested reader to the excellent monograph [Casson \& Bleilier 1988]  \cite{CB}. 
It is proved that each automorphism 
$\varphi\in Mod~(X)$ is isotopic to an automorphism $\varphi'$, such that either:
(i) $\varphi'$ has finite order; (ii) $\varphi'$ is an infinite order parabolic automorphism and preserves
 a Strebel foliation;   (iii) $\varphi'$ is of infinite order and does not preserve the simple closed curves
on $X$. In this case  there exist a stable ${\cal F}_s$
and unstable ${\cal F}_u$ mutually orthogonal measured foliations on $X$
such that $\varphi'({\cal F}_s)={1\over\lambda_{\varphi'}}{\cal F}_s$ 
and $\varphi'({\cal F}_u)=\lambda_{\varphi'}{\cal F}_u$, where $\lambda_{\varphi'}>1$
is called the {\it dilatation} of $\varphi'$. The foliations ${\cal F}_s,{\cal F}_u$ are minimal,
uniquely ergodic and describe $\varphi'$ up to a power.  
The automorphism $\varphi'$ is hyperbolic and invariant foliations ${\cal F}_s, {\cal F}_u$ are
 called {\it pseudo-Anosov foliations}.

\section{Proofs}
\subsection{Proof of lemma 1}
(i) For brevity, let ${\goth m}:=Jac~({\cal F})$;
suppose that $A=(a_{ij})\in GL_n({\Bbb Z})$ and 
$
\gamma_i'=\sum_{j=1}^na_{ij}\gamma_j
$
is a new basis in $H_1(X, Sing~{\cal F}; {\Bbb Z})$. We have to show, that ${\goth m}$
does not depend on the new basis; indeed,  using the integration rules we get:
\begin{equation}
 \lambda_i'  = \int_{\gamma_i'}\phi = \int_{\sum_{j=1}^na_{ij}\gamma_j}\phi=
\sum_{j=1}^n\int_{\gamma_j}\phi  = \sum_{j=1}^na_{ij}\lambda_j.
\end{equation}
To prove that ${\goth m}={\goth m}'$, consider the following equations:
\begin{equation}
{\goth m}'  = \sum_{i=1}^n{\Bbb Z}\lambda_i' = \sum_{i=1}^n {\Bbb Z} \sum_{j=1}^n a_{ij}\lambda_j=
\sum_{j=1}^n \left(\sum_{i=1}^n a_{ij}{\Bbb Z}\right)\lambda_j  \subseteq  {\goth m}. 
\end{equation}
Let $A^{-1}=(b_{ij})\in GL_n({\Bbb Z})$ be an inverse to the matrix $A$.
Then $\lambda_i=\sum_{j=1}^nb_{ij}\lambda_j'$ and 
\begin{equation}
{\goth m}  = \sum_{i=1}^n{\Bbb Z}\lambda_i = \sum_{i=1}^n {\Bbb Z} \sum_{j=1}^n b_{ij}\lambda_j'=
 \sum_{j=1}^n \left(\sum_{i=1}^n b_{ij}{\Bbb Z}\right)\lambda_j'  \subseteq  {\goth m}'.
\end{equation}
Since both ${\goth m}'\subseteq {\goth m}$ and ${\goth m}\subseteq {\goth m}'$, we conclude
that ${\goth m}' = {\goth m}$. Part (i) of lemma \ref{lm1} follows.

\medskip\noindent
(ii) Let $h: X\to X$ be an automorphism of surface $X$. Denote
by $h_*$ its action on $H_1(X, Sing~({\cal F}); {\Bbb Z})$
and by $h^*$ on $H^1(X; {\Bbb R})$ connected  by the formula: 
\begin{equation}
\int_{h_*(\gamma)}\phi=\int_{\gamma}h^*(\phi), \quad\forall\gamma\in H_1(X, Sing~({\cal F}); {\Bbb Z}), 
\quad\forall\phi\in H^1(X; {\Bbb R}).
\end{equation}
Let $\phi,\phi'\in H^1(X; {\Bbb R})$ be the closed forms, whose
trajectories define the foliations ${\cal F}$ and ${\cal F}'$, respectively.
Since ${\cal F}, {\cal F}'$ are topologically conjugate,
$
\phi'= \mu ~h^*(\phi)
$
for a $\mu>0$. 
Let $Jac~({\cal F})={\Bbb Z}\lambda_1+\dots+{\Bbb Z}\lambda_n$ and 
$Jac~({\cal F}')={\Bbb Z}\lambda_1'+\dots+{\Bbb Z}\lambda_n'$. Then:
\begin{equation}
\lambda_i'=\int_{\gamma_i}\phi'=\mu~\int_{\gamma_i}h^*(\phi)=
\mu~\int_{h_*(\gamma_i)}\phi, \qquad 1\le i\le n.
\end{equation}
By item (i), it holds:
\begin{equation}
Jac~({\cal F})=\sum_{i=1}^n{\Bbb Z}\int_{\gamma_i}\phi=
\sum_{i=1}^n{\Bbb Z}\int_{h_*(\gamma_i)}\phi.
\end{equation}
Therefore:
\begin{equation}
Jac~({\cal F}')=\sum_{i=1}^n{\Bbb Z}\int_{\gamma_i}\phi'=
\mu~\sum_{i=1}^n{\Bbb Z}\int_{h_*(\gamma_i)}\phi=\mu~Jac~({\cal F}).
\end{equation}
Lemma \ref{lm1} follows.
$\square$

\subsection{Proof of theorem 1}
We shall split the proof in a series of lemmas, starting with the following 
elementary
\begin{lem}\label{lm3}
Let $T\in M_n({\Bbb Z})$ be an endomorphism of the vector space ${\Bbb R}^n$,
where $M_n({\Bbb Z})$ is the set of the $n\times n$ matrices over ${\Bbb Z}$. 
If $Tx=\lambda x$, where $x\in {\Bbb R}^n$ is an eigenvector and $\lambda\in K$
is an eigenvalue of $T$, then $x$ can be scaled so that $x\in K^n$, where
$K={\Bbb Q}(\lambda)$ is a subfield of ${\Bbb R}$ generated by $\lambda$.  
\end{lem}
{\it Proof.} 
The proof is left  to the reader as an exercise in linear algebra;   otherwise,
we refer to the results of   [Borevich \& Shafarevich 1966,  Chapter 2, Section 2]  \cite{BS}.
$\square$

\bigskip\noindent
Let $f_N\in S_2(\Gamma_0(N))$ be a (normalized) Hecke eigenform, such that
$f_N(z)=\sum_{n=1}^{\infty}c_n(f_N)q^n$ its Fourier series. We shall denote
by $K_f={\Bbb Q}(\{c_n(f_N)\})$ the algebraic number field generated
by the Fourier coefficients of $f_N$. If  $g$ is the genus of the modular curve$X_0(N)$,
then $1\le deg~(K_f~|~{\Bbb Q})\le g$ and $K_f$ is a totally real field [Darmon 2004]  \cite{Dar1}, p. 25.
Let ${\cal F}_N$ be the foliation by vertical trajectories of $f_N$. The following lemma will be
critical. 
\begin{lem}\label{lm4}
$Jac~({\cal F}_N)$ is a ${\Bbb Z}$-module in the field $K_f$. 
\end{lem}
{\it Proof.} 
Let $\phi_N=\Re~(f_Ndz)$, where $f_N$ is the Hecke eigenform.  By the definition,
$Jac~({\cal F}_N)=\int_{H_1(X_0(N), Sing~{\cal F}_N; {\Bbb Z})}\phi_N=
{\Bbb Z}\lambda_1+\dots+{\Bbb Z}\lambda_g, ~\lambda_i\in {\Bbb R}$.
Let us study the action of the Hecke operators on the $Jac~({\cal F}_N)$. Since
$f_N$ is an eigenform, $T_n f_N=c_n f_N, ~\forall T_n\in {\Bbb T}_{\Bbb Z}$.
By virtue of the isomorphism $S_2(\Gamma_0(N))\cong \Omega_{hol}(X_0(N))$,
one gets $T_n\omega_N=c_n\omega_N$, where $c_n\in K_f$.
Let us evaluate the real parts of the last equation as follows:
$ \Re~(T_n\omega_N)=T_n(\Re~\omega_N)=\Re~(c_n\omega_N)=c_n(\Re~\omega_N)$;
note that the equality  $\Re~(c_n\omega_N)=c_n (\Re~\omega_N)$ involves the 
fact that $c_n$ is a real number.  One concludes that
$T_n\phi_N=c_n\phi_N, ~\forall T_n\in {\Bbb T}_{\Bbb Z}$ and  $c_n\in K_f$;
therefore,  the action of the Hecke operator $T_n$ on the $Jac~({\cal F}_N)$
can be written as:
\begin{equation}\label{eq7}
T_n(Jac~({\cal F}_N))=\int_{H_1}T_n\phi_N=\int_{H_1}c_n\phi_N=c_n Jac~({\cal F}_N),
\qquad c_n\in K_f,
\end{equation}
where $H_1=H_1(X_0(N), Sing~{\cal F}; {\Bbb Z})$ is the relative 
homology group. Thus, $T_n$ multiplies  $Jac ~({\cal F}_N)$ by 
the real number $c_n$.

Consider the ${\Bbb Z}$-module  $Jac~({\cal F}_N)={\Bbb Z}\lambda_1+\dots+{\Bbb Z}\lambda_g$, 
where $\lambda_i\in {\Bbb R}$.
The  Hecke operator $T_n$ is an endomorphism  of  the vector space 
$\{\lambda\in {\Bbb R}^g~|~\lambda=(\lambda_1,\dots,\lambda_g)\}$.  By virtue of 
(\ref{eq7}), vector $\lambda\in {\Bbb R}^g$ 
is an eigenvector of the linear operator $T_n\in M_g({\Bbb Z})$. 
In other words, $T_n\lambda=c_n\lambda, ~c_n\in K_f$.
In view of lemma \ref{lm3}, one can scale the 
vector $\lambda$ so that $\lambda_i\in K_f$.
$\square$

\begin{lem}\label{lm5}
$rank~(Jac~({\cal F}_N))= deg~(K_f~|~{\Bbb Q})$.
\end{lem}
{\it Proof.} 
Let $I_f=\{T\in {\Bbb T}_{\Bbb Z} : ~Tf=0\}$  and $R= {\Bbb T}_{\Bbb Z}/I_f$. 
It is well known that $R$ is an order in the number field $K_f$, i.e. a subring
of the ring of integers of $K_f$ containing $1$. Recall that a coefficient ring
of the lattice $\Lambda$ in  $K_f$ is defined as the set $\alpha\in K_f$, such that
$\alpha\Lambda\subseteq\Lambda$. The coefficient ring is always an order in $K_f$. 
It is verified directly, that $R$ is the coefficient ring of $Jac~({\cal F}_N)$.     
Moreover, $Jac~({\cal F}_N)$ can be scaled so that $Jac~({\cal F}_N)\subseteq R$
[Borevich \& Shafarevich 1966]  \cite{BS}, p.88. Therefore,  $Jac~({\cal F}_N)$ and $R$ have the same number of
generators, i.e. the same rank {\it ibid.}  But $rank~(R)= deg~(K_f~|~{\Bbb Q})$
[Diamond \& Shurman 2005]  \cite{DS}, p.234. Thus, $rank~(Jac~({\cal F}_N))= deg~(K_f~|~{\Bbb Q})$.
$\square$

\begin{lem}\label{lm6}
The foliation ${\cal F}_N$ is topologically conjugate to:

\medskip
(i) a Strebel foliation, if $deg~(K_f~|~{\Bbb Q})=1$;

\smallskip
(ii) a pseudo-Anosov foliation, if   $deg~(K_f~|~{\Bbb Q})=g$;

\smallskip
(iii) a degenerate pseudo-Anosov foliation, if $2\le deg~(K_f~|~{\Bbb Q})\le g-1$.
\end{lem}
{\it Proof.} We shall use the well known construction of measured 
foliations from the generators $\lambda_i>0$ of their Jacobians;
the construction is known as a {\it zippering of the rectangles}
[Veech 1982]  \cite{Vee1} and runs (roughly) as follows. 
Let $\lambda_1,\dots,\lambda_n$ be positive generators of $Jac~({\cal F})$
and fix a permutation $\pi$ of $n$ symbols; let $R_i$ be $n$ rectangles
with the base $\lambda_i$  placed one-next-to-other
in ${\Bbb R}^2$ and foliated by the vertical lines. One glues the
bottom $\lambda_i$ of $R_i$ with the top $\lambda_{\pi(i)}$ of 
$R_{\pi(i)}$; since the permutation $\pi$ preserves the total
length $\lambda_1+\dots+\lambda_n$ of the intervals,  all rectangles
will be glued in this way.  By zippering of the remaining vertical
sides of $R_i$ one gets a compact surface $X$ with the measured 
foliation ${\cal F}$ on it;  by the construction $Jac~({\cal F})$
is generated by $\lambda_i$.

In fact,  the Veech construction gives a little more than just a foliation on 
compact surface;   the surface comes with a complex structure,
see  [Veech 1982]  \cite{Vee1}, p. 215.   Our argument does not 
use the complex structure  but the reader can verify (given the arithmetic
nature of periods $\lambda_i$) that the underlying Riemann surface
is arithmetic, i.e. the orbit space of an arithmetic Fuchsian group.

\medskip
(i) Let  $deg~(K_f~|~{\Bbb Q})=1$, i.e the minimal possible. By lemma \ref{lm4}, 
\linebreak
 $rank~(Jac~({\cal F}_N))=1$
and therefore scaled   $\lambda_i$ are all  rational. Since $\lambda_i$ are
linearly dependent over ${\Bbb Q}$, the zippered rectangles will produce a 
measured foliation all of whose leaves, but a finite number, are closed. Such  a  foliation
is traditionally called the Strebel foliation [Strebel 1984]  \cite{S}.

\smallskip
(ii)  Let  $deg~(K_f~|~{\Bbb Q})=g$ be the maximal possible;
in this case $\lambda_i$ are linearly independent over ${\Bbb Q}$.
Thus, ${\cal F}_N$ is a  minimal foliation, i.e.
each non-singular leaf of ${\cal F}_N$  is dense on the surface 
$X_0(N)$. To show  that ${\cal F}$ is a pseudo-Anosov foliation,  
let $\lambda_i>0$ and $A\in GL_g({\Bbb Z})$ be a positive matrix
with the Perron-Frobenius eigenvector $(\lambda_1,\dots,\lambda_g)$. 
Denote by $i$ an involution, which acts on $S_2(\Gamma_0(N))$ according
with the formula
\begin{equation}
f(z)=\sum c_nq^n\longmapsto
f^*(z)=\sum \bar c_nq^n,
\end{equation}
where the bar sign means complex conjugation; since 
\linebreak
$S_2(\Gamma_0(N))\cong H^1(X_0(N); {\Bbb R})$, the 
involution defines a  split  
$H^1(X_0(N); {\Bbb R})= H^1_{\Bbb R}\oplus H^1_{i{\Bbb R}}$.
The matrix $A$ acts on  $H^1_{\Bbb R}$ and extends to the whole
space $H^1(X_0(N); {\Bbb R})$ by the involution $i$; the corresponding Perron-Frobenius 
eigenvector has the form $(\lambda_1,\dots,\lambda_g,\lambda_1,\dots,\lambda_g)$. 
By  the argument of [Thurston 1988]   
 \cite{Thu1} on p.427,   one concludes  that there
exists a pseudo-Anosov automorphism $\varphi_A: X_0(N)\to X_0(N)$,
such that $\varphi_A^*=A$; in other words, ${\cal F}_N$ is an
invariant foliation of the automorphism $\varphi_A$.

\smallskip
(iii) Let  $deg~(K_f~|~{\Bbb Q})=n$, where $2\le n\le g-1$;
in this case there are only $n$ rationally independent $\lambda_i$
among $\lambda_1,\dots,\lambda_g$. The foliation ${\cal F}_N$ is 
still minimal, since $n>1$;  however, it will acquire  $g-n$ extra separatrix
connections, i.e leaves joining the saddle points of ${\cal F}_N$. We  repeat the 
argument of item (ii) to prove that ${\cal F}_N$ is a pseudo-Anosov
foliation,  but  because of presence of the separatrix connections
it will be a degenerate one.           
$\square$

\medskip
Theorem \ref{thm1} follows from lemma \ref{lm6}. 
$\square$

\bigskip\noindent
{\sf Acknowledgments.}  
I am grateful to  Lawrence ~D. ~Taylor (University of  Nottingham)  for useful  correspondence
and to a referee for helpful comments.  



\vskip1cm

\textsc{The Fields Institute for Research in Mathematical Sciences, Toronto, ON, Canada,  
E-mail:} {\sf igor.v.nikolaev@gmail.com}

\smallskip
{\it Current address: Department of Mathematics, 
University of Sherbrooke, {\sf Igor.Nikolaev@USherbrooke.ca}}

\end{document}